\documentclass{article}

\usepackage{url}

\begin{document}

\title{Comment on ``Noether's-type theorems on time scales'' 
[J. Math. Phys. 61, 113502 (2020)]\thanks{This is a preprint 
of a paper whose final and definite form is published 
in the 'Journal of Mathematical Physics'.
Submitted 08-Jul-2022; Revised 10-Sep-2022; Accepted 13-Sep-2022.}}

\author{Delfim F. M. Torres\thanks{The author
is supported by FCT under project UIDB/04106/2020.}\\
{\tt delfim@ua.pt}}

\date{Center for Research and Development in Mathematics and Applications (CIDMA),
Department of Mathematics,\\ 
University of Aveiro, 3810-193 Aveiro, Portugal}

\maketitle

\begin{abstract}
We comment on the validity of Noether's theorem
and on the conclusions of [J. Math. Phys. 61 (2020), no.~11, 113502].
\end{abstract}

In \cite{[Anerot:Cresson:Belgacem:Pierret:2020]}, Anerot et al.\  
conclude that the Noether-type theorem proved by Bartosiewicz and Torres  
in \cite{[Bartosiewicz:Torres:2008]} is true along the 
Euler--Lagrange extremals under the assumption
that the Euler--Lagrange extremals are also satisfying the second Euler--Lagrange
equation: see \cite[bottom of page 9]{[Anerot:Cresson:Belgacem:Pierret:2020]}.
While this is true, the authors of \cite{[Anerot:Cresson:Belgacem:Pierret:2020]}
seem not aware that:
(i) such conclusion is not new (see \cite{[Torres:2004]});
(ii) in the class of $C^2$-functions where Noether's theorem 
is proved, all Euler--Lagrange extremals
do satisfy the second Euler--Lagrange equation (sometimes
also called the DuBois-Reymond condition): see also
\cite{[Torres:2004]} and references therein.

Indeed, it has been shown in \cite{[Torres:2004]} that while Noether's theorem is valid 
along every $C^2$ Euler--Lagrange extremal
(and every $C^2$ Euler--Lagrange extremal is also a DuBois-Reymond extremal),
Noether's result is not true when one enlarges the class
$C^2$ of admissible functions to be the class $Lip$ of Lipschitz functions,
for which the Euler--Lagrange equation is still valid.
Moreover, it has been proved in \cite{[Torres:2004]} that one can still
prove a version of Noether's theorem in the class of $Lip$ admissible functions,
but in that case we need to restrict the set of Euler--Lagrange
extremals to those that also satisfy the DuBois-Reymond condition
(the second Euler--Lagrange equation). So the conclusion 
of \cite{[Anerot:Cresson:Belgacem:Pierret:2020]}, that Noether's 
theorem is only true for Euler--Lagrange extremals
that also satisfy the second Euler Lagrange equation, is not 
a new result and does not invalidate, \emph{per se}, 
the Noether theorem proved in  \cite{[Bartosiewicz:Torres:2008]}: 
for the class of smooth admissible
functions considered by Noether's theorem, 
the second Euler--Lagrange equation is always satisfied 
along the Euler--Lagrange extremals.

However, the main message of \cite{[Anerot:Cresson:Belgacem:Pierret:2020]}
is that the result of \cite{[Bartosiewicz:Torres:2008]} is not correct.
To show that, an example is given. While the example is also not new, 
since it was borrowed from \cite{[Bartosiewicz:Torres:2008]}, we fully
agree that the example deserves indeed some comments to make things clear.
Unfortunately, the analysis of \cite{[Anerot:Cresson:Belgacem:Pierret:2020]} 
has some inconsistencies that do not help to clarify things: in the simulations
of \cite{[Anerot:Cresson:Belgacem:Pierret:2020]}, the authors
consider the initial condition $x(0) = x_0 = 0$
and a nonstandard condition $\Delta x_0 = 0.1$ that has never
been considered in the literature of calculus of variations
under Noether's theorem. 

Since paper \cite{[Anerot:Cresson:Belgacem:Pierret:2020]} is very technical, 
and not all readers may be familiar with the language of the calculus
of variations on time scales \cite{[Bohner:2004]}, 
to make things completely clear and understandable,
we restrict our discussion here to the classical case when 
the time scale is the set of real numbers. This allow us to easily explain 
the discrepancies between Noether's theorem and that particular example.
In our opinion such discrepancies are intrinsic to the nature of the 
calculus of variations, and do not depend on the theory of time scales.

The example of the calculus of variations given in \cite{[Bartosiewicz:Torres:2008]},
and recalled in \cite{[Anerot:Cresson:Belgacem:Pierret:2020]}, has a Lagrangian
given by $L(t,x,v) = \frac{x^2}{t} + t v^2$, defined 
in a time interval $t \in [a,b]$ that does not include $t = 0$. 
This means that for this Lagrangian the Euler--Lagrange equation
$$
\frac{d}{dt}\frac{\partial L}{\partial v}(t,x(t),x'(t)) 
= \frac{\partial L}{\partial x}(t,x(t),x'(t))
$$
reduces to
$$
\frac{d}{dt}\left(t x'(t)\right) = \frac{x(t)}{t},
$$
from which one concludes that the Euler--Lagrange extremals have the form
$$
x(t) = \frac{c_1}{t} + c_2 t,
$$
where the constants $c_1$ and $c_2$ are determined from 
the boundary conditions of the problem.
For example, if one considers the (nonstandard) boundary conditions $x(0) = 0$
and $x'(0) = 0.1$ considered in \cite{[Anerot:Cresson:Belgacem:Pierret:2020]},
then the Euler--Lagrange extremal is obtained with $c_1 = 0$ and $c_2 = 0.1$, that is,
the Euler--Lagrange extremal is the function $x(t) = 0.1 t$. However,
for standard boundary conditions $x(a) = \alpha$ and $x(b) = \beta$ of the calculus of variations
with $0 < a < b$, as considered in \cite{[Bartosiewicz:Torres:2008]},
then the Euler--Lagrange extremal is given by
$$
x(t) =  \frac{ab(b \alpha - a \beta)}{(b^2 - a^2)t} 
+ \frac{b \beta - a \alpha}{b^2 - a^2} t.
$$
What is important to mention here is that
$x(t)$ is an absolutely continuous function
that is not Lipschitz (and in particular does not belong to the class $C^2$ 
where Noether's theorem is proved). We also recall that while the 
Euler--Lagrange equation can still be proved in the class of Lipschitz functions, 
the Euler--Lagrange equation (and thus also Noether's theorem) is not valid 
in the class of absolutely continuous functions:
see, e.g., \cite{[Sarychev:Torres:2000]} and references therein.

\subsection*{Conflict of Interest}

The author has no conflicts to disclose.




\begin{thebibliography}{xx}

\bibitem{[Anerot:Cresson:Belgacem:Pierret:2020]}
B. Anerot, J. Cresson, K. Hariz Belgacem\ and\ F. Pierret,
Noether's-type theorems on time scales, 
J. Math. Phys. {\bf 61} (2020), no.~11, 113502, 31~pp. 
\url{https://doi.org/10.1063/1.5140201}

\bibitem{[Bartosiewicz:Torres:2008]}
Z. Bartosiewicz\ and\ D. F. M. Torres, 
Noether's theorem on time scales, 
J. Math. Anal. Appl. {\bf 342} (2008), no.~2, 1220--1226. 
\url{https://doi.org/10.1016/j.jmaa.2008.01.018}
{\tt arXiv:0709.0400}

\bibitem{[Bohner:2004]}
M. Bohner, 
Calculus of variations on time scales, 
Dynam. Systems Appl. {\bf 13} (2004), no.~3-4, 339--349. 

\bibitem{[Sarychev:Torres:2000]}
A. V. Sarychev\ and\ D. F. M. Torres, 
Lipschitzian regularity of minimizers for optimal 
control problems with control-affine dynamics, 
Appl. Math. Optim. {\bf 41} (2000), no.~2, 237--254.
\url{https://doi.org/10.1007/s002459911013}

\bibitem{[Torres:2004]}
D. F. M. Torres, 
Proper extensions of Noether's symmetry theorem 
for nonsmooth extremals of the calculus of variations, 
Commun. Pure Appl. Anal. {\bf 3} (2004), no.~3, 491--500. 
\url{https://doi.org/10.3934/cpaa.2004.3.491}

\end{thebibliography}
\end{document}